\documentclass{gtmon_a}
\pdfoutput=1

\usepackage{amscd}

\proceedingstitle{Groups, homotopy and configuration spaces (Tokyo
  2005)}
\conferencestart{5 July 2005}
\conferenceend{11 July 2005}
\conferencename{Groups, homotopy and configuration spaces,
                in honour of Fred Cohen's 60th birthday}
\conferencelocation{University of Tokyo, Japan}

\editor{Norio Iwase}
\givenname{Norio}
\surname{Iwase}

\editor{Toshitake Kohno}
\givenname{Toshitake}
\surname{Kohno}

\editor{Ran Levi}
\givenname{Ran}
\surname{Levi}

\editor{Dai Tamaki}
\givenname{Dai}
\surname{Tamaki}

\editor{Jie Wu}
\givenname{Jie}
\surname{Wu}

\title{Filtering the fiber of the pinch map}

\author{Brayton Gray}
\givenname{Brayton}
\surname{Gray}
\address{Dept of Math, Stats and Comp Sci (M/C 249)\\
University of Illinois at Chicago\\\newline
851 South Morgan Street\\Chicago, IL 60607-7045\\USA}
\email{brayton@uic.edu}

\volumenumber{13}
\issuenumber{}
\publicationyear{2008}
\papernumber{09}
\startpage{203}
\endpage{227}

\doi{}
\MR{}
\Zbl{}


\received{24 July 2006}
\revised{19 January 2007}
\accepted{26 January 2007}
\published{22 February 2008}
\publishedonline{22 Februray 2008}
\proposed{}
\seconded{}
\corresponding{}
\version{}

\arxivreference{0804.1896}

\makeatletter
\def\cnewtheorem#1[#2]#3{\newtheorem{#1}{#3}[section]
\expandafter\let\csname c@#1\endcsname\c@subsection}

\AtBeginDocument{\let\bar\wbar\let\tilde\wtilde\let\hat\what}

\makeop{im}

\swapnumbers
\cnewtheorem{theorem}[subsection]{Theorem}
\cnewtheorem{corollary}[subsection]{Corollary}
\cnewtheorem{lemma}[subsection]{Lemma}
\cnewtheorem{proposition}[subsection]{Proposition}
\theoremstyle{remark}
\cnewtheorem{definition}[subsection]{Definition}
  \let\c@equation\c@subsection
  \def\theequation{\bf \thesubsection}
\makeatother

\def\lraw{\longrightarrow}
\def\ome{\omega}
\def\Ome{\Omega}
\def\tms{\times}
\def\alp{\alpha}
\def\eps{\epsilon}
\def\wdg{\wedge}
\def\bta{\beta}
\def\sig{\sigma}
\def\Sig{\Sigma}
\def\gam{\gamma}
\def\dta{\delta}

\def\Dta{\Delta}
\def\wtit{\widetilde}
\def\pat{\partial}
\def\ol{\overline}
\def\noi{\noindent}

\def\Tht{\Theta}

\def\mbk{\medbreak}

\def\Raw{\Rightarrow}

\begin{document}

\begin{asciiabstract}

\end{asciiabstract}

\begin{htmlabstract}
A new type of Hopf invariant is described for the fiber of the pinch map:
<br>
F &rarr; X &cup;CA &rarr;<sup>&pi;</sup> SA
<br>
and this is used to study the boundary map in the fibration sequence
of Cohen, Moore and Neisendorfer:
<br>
&Omega;<sup>2</sup>S<sup>2n+1</sup>
&rarr;<sup>&part;<sub>n</sub></sup>
&Omega;F<sub>n</sub>  &rarr; &Omega;P<sup>2n+1</sup> &rarr; 
&Omega;S<sup>2n+1</sup>.
<br>
The boundary map is shown to be compatible with the Hopf invariant and a filtration of the spliting is obtained.
\end{htmlabstract}

\begin{abstract}
A new type of Hopf invariant is described for the fiber of the pinch map:
\[
F\longrightarrow X\cup CA \stackrel{\pi}{\longrightarrow} SA
\]
and this is used to study the boundary map in the fibration sequence
of Cohen, Moore and Neisendorfer:
\[
\Omega^2S^{2n+1} \stackrel{\partial_n}{\longrightarrow} \Omega F_n  \longrightarrow \Omega P^{2n+1} \longrightarrow \Omega S^{2n+1}.
\]
The boundary map is shown to be compatible with the Hopf invariant and a filtration of the spliting is obtained.
\end{abstract}

\maketitle

\section{Introduction}
\label{sec1}

In \cite{CMN1}, Cohen, Moore and Neisendorfer analyzed the homotopy
type of the mod $p^r$ Moore space
\[
P^{2n+1} = S^{2n} \cup_{p^r} e^{2n+1}
\]
by constructing a fibration sequence
\[
\Ome^2S^{2n+1} \stackrel{\pat_n}{\lraw} \Ome F_n  \lraw \Ome P^{2n+1} \lraw \Ome S^{2n+1}.
\]
A key result is that both $\Ome F_n$ and $\Ome P^{2n+1}$ contain a factor which is the loop space on a one point union of Moore spaces that are at least $4n-2$ connected.  Removing this factor (which can be inductively considered later) leaves a fibration sequence:
\[
\cdots \to \Ome^2 S^{2n+1} \stackrel{\ol\pat}{\lraw} \ol F_n \lraw T^{2n+1} \lraw \Ome S^{2n+1}
\]
and they show that 
\[
\ol F_n \cong S^{2n-1} \tms \prod\limits_{k\ge 1} S^{2n p^k -1} \{ p^{r+1}\}.
\]
An important remaining question is to understand the components of $\ol\pat$:
\begin{eqnarray*}
\ol\pat_0\co \Ome^2 S^{2n+1} &\lraw & S^{2n-1}\\
\ol\pat_k\co \Ome^2 S^{2n+1} &\lraw & S^{2np^k -1}\{ p^{r+1}\}.
\end{eqnarray*}
$\ol\pat_0$ is well understood and plays a key role in homotopy theory.  
This is the noted ``Cohen--Moore--Neisendorfer map''.  It fits into a homotopy commutative diagram:
\[
\begin{array}{ccc}
S^{2n-1}      &\stackrel{p}{\lraw}      &S^{2n-1} \\
\bigg\downarrow &{}^{\scriptstyle\ol\pat_0}{\Huge{\nearrow}}  &\bigg\downarrow   \\
\Ome^2 S^{2n+1}  &\stackrel{p}{\lraw}   &\Ome^2 S^{2n+1} \\
\end{array}
\]
the consequence of which is that
\[
p \cdot \pi_k (S^{2n+1}) \subset E^2 \pi_{k-2} (S^{2n-1}) \qquad k>0
\]
from which the exponent result follows:
\[
p^n \pi_k (S^{2n+1})=0 \qquad k>0.
\]
There is much interest in understanding $\ol\pat_k$ for $k>0$; in
particular, it is not known whether they are all null homotopic (see
Anick and Gray \cite{AG}, Cohen, Moore and Neisendorfer \cite{CMN1},
Gray \cite{G2,G3} and Gray and Theriault \cite{GT}).  In case $r>1$,
Neisendorfer \cite{N1,N2} has shown that $\ol\pat_k=0$ for
$k\ge 1$, so our emphasis will be on the case $r=1$.

A remarkable observation of Cohen, Moore and Neisendorfer \cite{CMN2}
is that there is an isomorphism of Hopf algebras (see \fullref{secAppendix}):
\[
H_* (\Ome^2 S^{2n+1} ; \mathbb{Z}/p)\cong H_* (\ol F_n; \mathbb{Z}/p).
\]
However, the induced homomorphism $(\ol\pat)_*=0$.

The intention of this work is to consider constructions involving $\ol
F_n$ which are analogous to constructions involving $\Ome^2S^{2n+1}$.
We will be able to do this compatibly with the map $\ol\pat$.  In
particular, we will look at the filtration of Selick \cite{S} and the
classifying space of the double suspension \cite{G2}.  Our results
will follow from the construction of a new ``Hopf invariant'' type
map:
\[
\begin{CD}
F_n @> {h} >> F_{np}
\end{CD}
\]
which induces an epimorphism in $p$--local homology.  Recall that there is a $p$--local fibration sequence:
\[
\begin{CD}
S^{2n}_{(p^s-1)} @>>> S^{2n}_\infty @> {H_{p^s}}>> S^{2np^s}_\infty.
\end{CD}
\]
Here we are using James' notation $X_\infty$ for the reduced product
space and $X_{(r)}$ for the subspace of words of length less than or
equal to $r$.  We also write $F_{(r)}$ for the analogous filtration of
the relative James construction for $F$ (Gray \cite{G1}).  By analogy,
we compare $F_{(p^s-1)}$ to the fiber $\Dta_s$ of the Hopf invariant
$h^s$:
\[
\begin{CD}
\Dta_s @>>> F_n @> {h^s}>> F_{np^s}
\end{CD}
\]
\begin{theorem}
$\Ome \Dta_s \simeq \Ome F_{(p^s-1)} \tms \Ome X_s$ where $X_s$ is a wedge of mod $p^r$ Moore spaces.
\label{thm1.1}
\end{theorem}

This will appear as \fullref{thm8.1}. These same ideas lead to:

\begin{theorem}
If $k>1$, there is a homotopy commutative diagram:
\[
\begin{CD}
\Ome^2 S^{2n+1} @> \ol\pat_{n(k)} >> S^{2np^k-1}\{p^{r+1}\} \\
 @V \gam_{_k} V \cong V      @AA \ol\pat_{np^{k-1}}(1) A \\
\Ome^2 S^{2n+1} @> \Ome H_{p^{k-1}} >> \Ome^2S^{2np^{k-1}+1}
\end{CD}
\]
where $H_{p^k}$ is a James--Hopf invariant and $\gam_{_k}$ is a homotopy equivalence.
\label{thm1.2}
\end{theorem}

\begin{theorem}
The map $\ol\pat_{n(k)}\co \Ome^2S^{2n+1}\lraw S^{2np^k-1} \{p^{r+1}\}$ is homotopic to a composition:
\[
\Ome^2 S^{2n+1} \stackrel{\nu}{\lraw} BW_n \stackrel{\epsilon}{\lraw} S^{2np^k-1}\{p^{r+1}\}.
\]
\label{thm1.3}
\end{theorem}

These results follow from equation \eqref{eq9.2} and \fullref{cor10.3}
respectively, in the sequel.

Throughout this paper all spaces will be assumed to be localized at a fixed prime $p>2$ and all homology will be $p$--local unless otherwise indicated.

\section{Filtration of $\Ome^2S^{2n+1}$}
\label{sec2}

$H_*(\Ome^2S^{2n+1}; \mathbb{Z}/p)$ is a free commutative algebra on
generators of dimensions $2np^i-1$, $2np^{i+1}-2$ for $i\ge 0$.
Selick \cite{S} has described a sequence of $H$--spaces whose homology
filters $H_*(\Ome^2S^{2n+1}; \mathbb{Z}/p)$ by successively adding one
generator at each stage.  Let us write $S^{2n}_{(k)}=J_k(S^{2n})$ for
the $k$th filtration of the James construction
$J(S^{2n})=S^{2n}_\infty$.  Then the space capturing the first $2s$
generators of $H_*(\Ome^2S^{2n+1}; \mathbb{Z}/p)$ is precisely $\Ome
S^{2n}_{(p^s-1)}$.  Selick \cite{S} defines spaces $G_s$ which contain
the first $2s+1$ generators.  We define $G_s$ by a diagram of
fibrations:

{\small\[
\begin{CD}
G_s  @> >>  S^{2np^s-1} @> \alp >> S^{2n}_{(p^s-1)} \\
@ V VV  @VVV  @| \\
\Ome^2 S^{2n+1} @> >> \Ome^2 S^{2np^s+1} @> >> S^{2n}_{(p^s-1)} @>> > \Ome S^{2n+1} @> H_{np^s} >> \Ome S^{2np^s+1}\\
@V VV  @VV \nu V \\
BW_{np^s} @= BW_{np^s} 
\end{CD}
\]}
where $\alp\co S^{2np^s-1}\to S^{2n}_{(p^s-1)}$ is the attaching map for
the $2np^s$ cell of $S^{2n}_\infty$, $H_{np^s}$ is the appropriate
James--Hopf invariant, and $BW_{np^s}$ is the classifying space for the
double suspension.  Since $G_s$ is the fiber of $\nu\circ \Ome
H_{np^s}$, it is an $H$ space if $p>3$ (Gray
\cite[Proposition~6]{G2}).  $G_s$ could also be described as the fiber
of a Toda--Hopf invariant $\Ome S^{2n}_{(p^{s+1}-1)}
\stackrel{H'}{\lraw} \Ome S^{2np^{s+1}-1}$ using the techniques of
Gray \cite{G2} and Moore and Neisendorfer \cite{MN}).

Constructing an analogous filtration for $\ol F_n$ is not difficult.  The key result is in the compatibility of the two filtrations.  

\section{Proof of Theorem CMN}
\label{sec3}

In this section we give a brief summary of the proof of the main
result of Cohen, Moore and Neisendorfer \cite{CMN1}.

\mbk\noi
{\bf Theorem CMN}\qua {\sl There is a diagram of fibration sequences:
\[
\begin{CD}
\Ome F_n @>>> \Ome P^{2n+1} @>>> \Ome S^{2n+1}\\
@VVV  @VVV  @|  \\
S^{2n-1}\tms V_n @>>> T^{2n+1} @>>> \Ome S^{2n+1}\\
@VV i V   @VV i' V \\
P_n  @> = >>   P_n \\
@VV \ol \phi V      @VVV \\
F_n @>>> P^{2n+1} @> p >> S^{2n+1}
\end{CD}
\]
where the maps $i$ and $i'$ are null homotopic, $p$ is the pinch map, $P_n$ is a one point union of mod $p^r$ Moore spaces of dimension $\ge 4n$ and
\[
V_n=\prod_{k\ge 1} S^{2np^k-1} \{p^{r+1}\}.
\]
Here $S^m\{d\}$ is the fiber of the degree $d$ map on $S^m$.  The inessentiality of $i$ and $i'$ implies that\/}
\begin{eqnarray*}
&&\Ome F_n\simeq S^{2n-1} \tms V_n\tms \Ome P_n\\
&&\Ome P^{2n+1} \simeq T^{2n+1} \tms \Ome P_n.
\end{eqnarray*}

\proof[Sketch of the proof]  Both the mod $p^r$ homotopy and the mod $p$ homology of $\Ome P^{2n+1}$ have differential Lie algebra structures, and the Hurewicz map is a Lie algebra homomorphism.  Furthermore $\Ome F_n$ has an extended ideal structure.  These structures are obtained from the Samelson product and the $r$th Bockstein $\bta^{(r)} (\bta^{(i)}=0 \ \mbox{\rm for } i<r)$.  

$H_*(\Ome P^{2n+1};\mathbb{Z}/p)=U(L)$ where $L$ is a free Lie algebra on $v\in H_{2n}(\Ome P^{2n+1}; \mathbb{Z}/p)$ and $u=\bta^{(r)}(v)$ and $U(L)$ is the universal enveloping algebra.  $L^{(0)}\subset L$ is the Lie ideal generated by $x_i=ad^{i-1} (v)(u)$ for $i\ge 1$ and $H_*(\Ome F_n;\mathbb{Z}/p)=U(L^{(0)})$.

Furthermore the suspensions  of the $x_i$,
\[
\sigma (x_i)\in H_{2ni} (F_n;\mathbb{Z}/p) \cong \mathbb{Z}/p
\]
are non zero.  $u$ and $v$ lie in the image of the mod $p^r$ Hurewicz homomorphism, so all of $L$ and hence $L^{(0)}$ does as well
\[
L^{(0)}\subset \im \{
 \pi_* (\Ome F_n; \mathbb{Z}/p^r)\lraw H_* (\Ome F_n; \mathbb{Z}/p)\}
\]
using the extended ideal structure in homotopy.

Both in homotopy and homology, $\bta^{(r)} (x_{p^k})=0$, so one may construct an extension:
\[
P^{2np^k-1}(p^{r+1})\stackrel{\dta_k}{\lraw}\Ome F_n
\]
of the mod $p^r$ homotopy class
\[
P^{2n p^k-1}\lraw \Ome F_n
\]
representing $x_{p^k}$ in homology.  The only property of the maps $\dta_k$ that is needed is that their Hurewicz image is $x_{p^k}$.  From these maps the authors construct
\[
S^{2np^k-1} \{p^{r+1}\}\stackrel{e}{\lraw} \Ome P^{2np^k} (p^{r+1}) \stackrel{\Ome\ol\dta_k}{\lraw}\Ome F_n
\] 
where $\ol\dta_k$ is the adjoint of $\dta_k$ and $e$ is a particular map described in the next section.  These maps are assembled via loop multiplication to obtain:
\[
\Tht\co  S^{2n-1} \tms V_n \lraw \Ome F_n
\]
The Bockstein $\bta^{(r)}$ is trivial in $H_* (S^{2n-1}\tms V_n;\mathbb{Z}/p)$ but not in $H_*(\Ome F_n;\mathbb{Z}/p)$.  However $\Tht$ induces an isomorphism in the Bockstein homology of these homology groups (see \fullref{secAppendix} for a calculation of the Bockstein homology).

Next the authors construct a sequence of sub-Lie algebras $L^{(k+1)}\subset L^{(k)}\subset L^{(0)}$ via short exact sequences of Lie algebras:
\[
\begin{array}{cccccccc}
0   &\lraw  &L^{(1)}      &\lraw  &L^{(0)}  &\lraw  &\langle \tau_0\rangle                   &\lraw 0\\
0   &\lraw  &L^{(k+1)}  &\lraw  &L^{(k)}  &\lraw  &\langle \tau_k, \sigma_k\rangle  &\lraw 0.
\end{array}
\]
Here $\langle \tau_0\rangle$ and $\langle \tau_{k_1} \sig_k\rangle$ are free commutative Lie algebras generated by $\tau_k=x_{p^k}$ of dimension $2np^k-1$ and
\[
\sig_k=\frac{1}{2p} \sum\limits^{p^k-1}_{i=1} \binom{p^k}{i} [x_i, x_{p^k -i}]\eps L^{(0)}.
\]
This is possible since $\tau_k$ is a generator and $\sig_k$ is decomposable in $L^{(k-1)}$ but not in $L^{(k)}$ since $\frac{1}{2p} \binom{p^k}{p^{k-1}}$ is a $p$--local unit and $x_{p^{k-1}}=\tau_{k-1}\notin L^{(k)}$.  

It follows that there is a split short exact sequence of differential Hopf algebras:
\[
0\lraw U (L^{(k+1)}) \lraw U(L^{(k)})\lraw H_* (S^{2np^k-1}\{p^{r+1}\};\mathbb{Z}/p)\to 0
\]
and hence 
\[
U(L^{(\infty)})\otimes H_p (S^{2n-1}\tms V_n;\mathbb{Z}/p)\cong U(L^{(0)})=H_* (\Ome F_n;\mathbb{Z}/p)
\]
where $L^{(\infty)}=\bigcap\limits_{k\ge 0} L^{(k)}$.

Consequently $H_t(U(L^{(\infty)});\bta^{(r)})=0$ for $t>0$ and hence $L^{(\infty)}$ has a basis consisting of classes $\{ x_\alp, \bta^{(r)} x_\alp\}$.  Since $x_\alp\in L^{(\infty)}\subset L^{(0)}$ one can choose maps
\[
\phi_\alp\co P^{n_\alp} \lraw \Ome F_n
\]
whose homology image is $x_\alp$ and $\bta^{(r)} x_\alp$ using the extended ideal structure.  Assembling these one produces a one point union of mod $p^r$ Moore spaces, $P_n$ and a map:
\[
\ol\phi\co P_n\lraw F_n
\]
such that the homology image of 
\[
(\Ome \ol\phi)_* \co  H_* (\Ome P_n;\mathbb{Z}/p)\lraw H_* (\Ome F_n;\mathbb{Z}/p)
\]
is exactly $U(L^{(\infty)})$. From this the authors produce a homotopy equivalence
\[
S^{2n-1} \tms V_n\tms\Ome P_n\lraw \Ome F_n
\]
by multiplying the maps $\Tht$ and $\ol\phi$ via the $H$ space structure in $\Ome F_n$.  Let $\ol F_n$ be the fiber of $\ol\phi$.  Then
\[
S^{2n-1}\tms V_n \stackrel{\Tht}{\lraw} \Ome F_n \stackrel{\pat}{\lraw} \ol F_n
\]
is a homotopy equivalence and $i$ is inessential.  This completes the outline of the proof.\endproof

\section{Combinatorial description of $F_n$}
\label{sec4}

In \cite{G1} a combinatorial description of the fiber of the pinch map:
\[
F\lraw X\cup CA \stackrel{\pi}{\lraw} SA
\]
was given in the spirit of the James construction $X_\infty$ for $\Ome
SX$.  The model, designated $(X,A)_\infty$ consists of all words in
$X_\infty$ with the property that all letters after the first letter
are required to lie in $A$, where $A\subset X$.  Alternatively, this
can be described by a push out diagram:
\begin{equation}
\begin{CD}
X\tms A_\infty @ >>> (X,A)_\infty\\
@AAA    @AAA\\
A\tms A_\infty @>>> A_\infty
\end{CD}
\label{eq4.1}
\end{equation}
\begin{proposition}
{\rm \cite{G1}}\qua There is a map $(X,A)_\infty \to F$ which is a homotopy equivalence when the inclusion $A\subset X$ is a cofibration. 
\label{prop4.2}
\end{proposition}
There is an action of the monoid $A_\infty$ on $(X,A)_\infty$ and $(X,A)_\infty$ can be thought of as a universal space in the following sense.  If $Y$ is any space on which $A_\infty$ acts and $g\co X\to Y$ is any map such that $g(a)=a\cdot *$ for some point $*\in Y$, there is a unique $A_\infty$ equivariant map
\[
g_\infty \co (X,A)_\infty \to Y.
\]
(See \cite[4.2]{G1}).  The map $(X,A)_\infty \to F$ is constructed from the action $A_\infty \tms F\to \Ome SA\times F\to F$.

The orbit space of $(X,A)_\infty$ under the action of $A_\infty$ is $X/A$ and we may use the universal property to establish the following diagram
\begin{equation}
\begin{CD}
(X,A)_\infty @> e >> X_\infty \\
@V \rho VV    @VVV \\
X/A @>>> (X/A)_\infty
\end{CD}
\label{eq4.3}
\end{equation}
Note that the inclusion $e\co (X/A)_\infty\to X_\infty$ is the unique $A_\infty$ equivariant extensions of the inclusion of $X$ in $X_\infty$.

\mbk\noi
{\bf Example 4.4}\qua  Let $A=S^{2n-2}\subset P^{2n-1}(p^t)=X$.  Then $(X,A)_\infty$ is the homotopy fiber of the map of degree $p^t$
\[
(X,A)_\infty = S^{2n-1}\{p^t\}\lraw S^{2n-1}\stackrel{p^t}{\lraw} S^{2n-1}
\]
and the map $e\co S^{2n-1}\{p^t\}\lraw (P^{2n-1} (p^t))_\infty \cong \Ome P^{2n} (p^t)$ is uniquely determined as a $\Ome S^{2n-1}$ equivariant map extending the inclusion of $P^{2n-1}(p^t)$.

One of the main features of the construction $(X,A)_\infty$ is that we can define functorial Hopf invariants using the same formulas as in James \cite{J}.  In particular we have a pointwise commutative diagram:
\[
\begin{CD}
A_\infty @ > H_k >> (A^{(k)})_\infty \\
@VVV   @VVV \\
(X,A)_\infty @> H_k >> (X\wedge A^{(k-1)})_\infty \\
@VVV   @VVV \\
X_\infty @> H_k >> (X^{(k)})_\infty 
\end{CD}
\]
It would be desirable to construct functorial compressions of these maps:
\[
h\co (X,A)_\infty \lraw (X\wedge A^{(k-1)}, A^{(k)})_\infty\subset (X\wedge A^{(k-1)})_\infty
\]
but I have been unable to do this.  In the next section we will construct a map $h$ of this form, but we have no knowledge of how it relates to $H_k$.  

\section{Construction of $h\co F_n\lraw F_{np}$}
\label{sec5}

In this section we will define a kind of Hopf invariant which is key for the results of this paper.

\begin{theorem}
Suppose $(X,A)$ is a suspension pair and $r\ge 1$.  Then there is a map:
\[
h\co (X,A)_\infty \lraw (X\wedge A^{(r-1)}, A^{(r)})_\infty
\]
such that the composition:
\[
X\tms A_\infty \lraw (X,A)_\infty \stackrel{h}{\lraw} (X\wedge A^{(r-1)}, A^{(r)})_\infty \stackrel{\rho}{\lraw} X\wedge A^{(r-1)}/A^{(r)}
\]
is homotopic to the composition: 
\[
X\tms A_\infty \stackrel{1\tms H_{r-1}}{\lraw}  X\tms A^{(r-1)}_\infty \lraw X\wedge A^{(r-1)}_\infty \stackrel{\alpha}{\lraw}  X\wedge A^{(r-1)} \lraw X\wedge A^{(r-1)}/A^{(r)}
\]
where $H_{r-1}\co A_\infty \to A^{(r-1)}_\infty$ is any map and $\alpha$ uses the suspension structure of $X$ to collapse $X\wedge A^{(r-1)}_\infty$ to $X\wedge A^{(r-1)}$.
\label{thm5.1}
\end{theorem}
\begin{corollary}
Let $F_n$ be the fiber of the pinch map
\[
p\co  P^{2n+1} (p^r)\lraw S^{2n+1}
\]
for each $n$ and $H_{r-1}$ be any choice of James--Hopf invariants.  Then there is a map:
\[
h\co F_n\lraw F_{rn}
\]
such that the composition: 
\[
S^{2n}\tms \Ome S^{2n+1} \lraw F_n \stackrel{h}{\lraw} F_{rn} \lraw P^{2rn+1}
\]
is homotopic to:
\[
S^{2n} \tms \Ome S^{2n+1} \stackrel{1\tms H_{r-1}}{\longrightarrow} S^{2n}\tms \Ome S^{2(r-1)n+1} \lraw S^{2rn} \lraw  P^{2rn+1}; 
\]
in particular, $h_*\co H_{2rn}(F_n)\lraw H_{2rn} (F_{rn})$ is an isomorphism.
\label{cor5.2}
\end{corollary}

{\bf Note}\qua It is an easy calculation to see that
\[
H_i(F_n)= \left\{
\begin{array}{cc}
\mathbb{Z}_{(p)} & n\mid i \\
0 &n\nmid i
\end{array}
\right. 
\]
and the map $F_n \stackrel{\pi}{\lraw} P^{2n+1}(p^r)$ is reduction mod $p^r$ in homology.  In particular we have defined a map:
\[
F_n \lraw P^{2nr+1} (p^r)
\]
for each $r\ge 1$ which is onto in $p$--local homology.

\mbk\noi
{\bf Proof of \fullref{thm5.1}}\qua  We begin by constructing a map 
\[
\mu\co (X,A)_\infty \lraw (X,A)_\infty/A_\infty \lraw X\tms A_\infty/A\rtimes A_\infty \equiv X\rtimes A_\infty/A\rtimes A_\infty
\]
which follows from the push out diagram (4.0).  This map is functorial and we have a commutative diagram:
\[
\begin{CD}
(X,A)_\infty @> \mu>> X\rtimes A_\infty/A\rtimes A_\infty \\
@AAA   @AAA \\
X\tms A_\infty @>>> X\rtimes A_\infty
\end{CD}
\]
using the functorial property, we have the commutative diagram:
\[
\begin{CD}
(X,A)_\infty @>>> X\rtimes A_\infty/A\rtimes A_\infty \\
@VVV   @VVV \\
(CX,A)_\infty @>>> CX\rtimes A_\infty/A\rtimes A_\infty \simeq S(A\rtimes A_\infty)
\end{CD}
\]
however $(CX,A)_\infty$ is the fiber of $CX/A\lraw SA$ and hence is contractible.  
Consequently there is a lifting $\gam$ of $\mu$:
\[  
\begin{array}{ccc}
(X\rtimes A_\infty , A\rtimes A_\infty)_\infty &\stackrel{\epsilon}{\lraw} &(X\wdg A_\infty, A\wdg A_\infty)_\infty\\
{}^{\textstyle \gam} {\textstyle \nearrow} \hfil\hfil\hfil \big\downarrow {\scriptstyle\rho}_{_1} \hfil\hfil\hfil & & \big\downarrow {\scriptstyle\rho}_{_2} \\
(X,A)_\infty \stackrel{\mu}{\lraw} X\rtimes A_\infty /A\rtimes A_\infty &\lraw &X\wdg A_\infty / A\wdg A_\infty.
\end{array}
\]
We now apply the homotopy commutative diagram:
{\small\[
\begin{CD}
(X\wdg A_\infty, A\wdg A_\infty)_\infty @> \bta >> (X\wdg A^{(r-1)}_\infty, A\wdg A^{(r-1)}_\infty)_\infty @ >\alp_\infty >> (X\wdg A^{(r-1)}, A^{(r)})_\infty\\
@VV \rho_{_2} V  @VV \rho_{_3} V  @VV \rho V\\
X\wdg A_\infty / A\wdg A_\infty @ > \bta' >> X\wdg A^{(r-1)}_\infty / A \wdg A^{(r-1)}_\infty @> \alp' >> X\wdg A^{(r-1)}/ A^{(r)}
\end{CD}
\]}
where $\bta$ and $\bta'$ are defined using $H_{r-1}\co  A_\infty \to A^{(r-1)}_\infty$.  We now define $h=\alp_\infty\bta\epsilon\gam$.  Clearly the composition:
\[
\begin{CD}
X\tms A_\infty @> >> (X,A)_\infty @> h >> (X\wdg A^{(r-1)}, A^{(r)})_\infty @> \rho >> X\wdg A^{(r-1)}/A^{(r)}
\end{CD}
\]
is homotopic to the composition on the left and top in the diagram:
\[
\begin{array}{cclcc}
X\wdg A_\infty/A\wdg A_\infty &\stackrel{\bta'}{\lraw} &X\wdg A^{(r-1)}_\infty /A\wdg A^{(r-1)}_\infty &\stackrel{\alp'}{\lraw} &X\wdg A^{(r-1)}/A^{(r)}\\
\bigg\uparrow  &&\qquad\bigg\uparrow  &&\bigg\uparrow \\
X\tms A_\infty &\stackrel{\scriptsize 1\tms H_{r-1}}{\displaystyle\lraw} &X\tms A^{(r-1)}_\infty \stackrel{\eps'}{\lraw}X\wdg A^{(r-1)}_\infty &\stackrel{\alp}{\lraw} &X\wdg A^{(r-1)}
\end{array}
\]
and consequently to the composite along the bottom and the right as well.\qed

\begin{proposition}
The map $h\co F_n\lraw F_{np^s}$ induces an epimorphism in $p$--local homology.
\label{prop5.3}
\end{proposition}

\mbk\noi
{\bf Proof}\qua  We use the map $\Ome S^{2n+1} \stackrel{\partial}{\lraw} F_n$ which has degree $p$ in $H^{2nk}$ for each $k>0$ to calculate the cup product structure in $H^*(F_n)$.  We choose generators $e_i\in H^{2ni} (F_n)$ such that $\pat^*(e_i)$ is $p$ times the generator of $H^{2ni}(\Ome S^{2n+1})$ dual to the $i$th power of a homology generator in $H_{2n}(\Ome S^{2n+1})$.  Then it is easy to see that
\[
e_ie_j=p\binom{i+j}{i} e_{i+j}.
\]
Let us designate $d_i\in H^{2nip^s}(F_{np^s})$ for the corresponding generator; then
\[
h^* (d_1) = u_1 e_{p^s}
\]
for some $p$ local unit $u_1$ by \fullref{cor5.2}.  We show that
\[
h^*(d_i)=u_i e_{ip^s}
\]
where $u_i$ is a $p$--local unit for each $i\ge 1$ by induction.  Using the product structure we have
\[
p id_i = d_1 d_{i-1},
\]
so $pi h^*(d_i)=h^*(d_1d_{i-1})=h^*(d_1)h^* (d_{i-1})$
\begin{eqnarray*}
&=& (u_1e_{p^s}) (u_{i-1}e_{(i-1)p^s})\\
&=& p u_1 u_{i-1} \binom{ip^s}{p^s} e_{ip^s}.
\end{eqnarray*}
It suffices to show that $\frac{1}{i} \binom{ip^s}{p^s}$ is a $p$--local unit.  

Now let $v_p(m)$ be the exponent of $p$ in $m$ and $[x]$ be the greatest integer less than or equal to $x$.  Then
\[
v_p (n!)=\sum\limits_{i\ge 1} \left[\frac{n}{p^i}\right].
\]
Consequently $v_p((p^si)!)=p^{s-1}i+p^{s-2}i+\cdots + i+v_p(i!)$ so 
\begin{eqnarray*}
v_p \left(\binom{p^si}{p^s}\right) &=& v_p((p^si)!) -v_p ((p^s(i-1))!)-v_p(p^s!)\\
&=& v_p(i!)-v_p((i-1)!)\\
&=& v_p(i).
\end{eqnarray*}
and we are done.\qed

By \fullref{cor5.2}, the composition
\[
\Ome S^{2n+1}\stackrel{\pat}{\lraw} F_n \stackrel{h}{\lraw} F_{np^s}\lraw P^{2np^2+1}
\]
is null homotopic, so there is a lifting $\wtit H_{p^s}$
\[
\begin{CD}
\Ome S^{2n+1} @ >\wtit H_{p^s} >> \Ome S^{2np^s+1} \\
@V \pat VV  @VV \pat V \\
F_n @ > h >> F_{np^s}.
\end{CD}
\]
\begin{proposition}
$(\wtit H_{p^s})_*\co H_* (\Ome S^{2n+1})\to H_* (\Ome S^{2np^s+1})$ is an epimorphism.
\label{prop5.4}
\end{proposition}

This follows since both maps labelled $\pat$ have
degree $p$ in homology and $h_*$ is an
epimorphism.\endproof

In particular, the fiber of $\wtit H_{p^s}$ is $S^{2n}_{(p^s-1)}$ as if $\wtit H_{p^s}$ were the James--Hopf invariant $H_{p^s}$.  Both $\Ome\wtit H_{p^s}$ and $\Ome H_{p^s}$ can be placed in the fibration sequence induced by the inclusion $S^{2n}_{(p^s-1)}\subset \Ome S^{2n+1}$, so there is an equivalence $\gam_s\co \Ome^2S^{2n+1}\to \Ome^2S^{2n+1}$ such that
\[
\Ome \wtit H_{p^s}=\Ome H_{p^s}\circ\gam_s
\]
hence we have:

\begin{corollary}
There is a homotopy commutative diagram
\[
\begin{CD}
\Ome^2 S^{2n+1} @> \pat_n >> \Ome F_n \\
@V \Ome \wtit H_{p^s} VV         @ VV \Ome h^s V \\
\Ome^2 S^{2np^s+1} @> \pat_{np} >> \Ome F_{np^s}
\end{CD}
\]
where $\Ome \wtit H_{p^s}\sim \Ome H_{p^s}\circ \gam_s$ for some homotopy equivalence $\gam_s\co \Ome^2 S^{2n+1}\lraw$\break $\Ome^2 S^{2n+1}$.  
\label{cor5.5}
\end{corollary}

\section{A Filtered Decomposition}
\label{sec6}

Let $\Dta_s$ be the homotopy fiber of $h^s\co F_n\to F_{np^s}$.  In this section we will compare the decompositions of $\Ome F_n $ and $\Ome F_{np^s}$ and prove

\begin{theorem}
$\Ome \Dta_s\cong S^{2n-1} \tms \prod\limits_{1\le k <s} S^{2np^k-1} \{p^{r+1}\}\tms\Ome S^{2np^s-1}\tms\Ome R_s$ where $R_s$ is a wedge of mod $p^r$ Moore spaces.
\label{thm6.1}
\end{theorem}

Let $\Psi\co S^{2n-1}\tms V_n\to S^{2np-1}\tms V_{np}$ be defined by first projecting onto $V_n$ and then applying  

\[
V_n = S^{2np-1}\{p^2\}\tms V_{np}\stackrel{\rho\tms 1}{\lraw} S^{2np-1} \tms V_{np}
\]
where $\rho\co  S^{2np-1}\{p^2\}\to S^{2np-1}$ is the projection.

\begin{lemma}
With appropriate choices of the maps $\Tht_n$ and $\Tht_{np}$ there is a homotopy commutative diagram:
\[
\begin{CD}
S^{2n-1} \tms V_n @> \Tht_n >> \Ome F_n \\
@VV\psi V  @VV \Ome h V \\
S^{2np-1} \tms V_{np} @> \Tht_{np} >> \Ome F_{np}.
\end{CD}
\]
\label{lem6.2}
\end{lemma}

{\bf Proof}\qua  Given choices $\ol\dta_k\co  P^{2np^k}(p^{r+1})\to F_n$ inducing epimorphisms in $p$--local cohomology, define $\ol\dta_k (np)$ as the composition:
\[
P^{2np^{k+1}}(p^{r+1}) \stackrel{\ol\dta_{k+1}}{\lraw} F_n \stackrel{h}{\lraw} F_{np}
\]
From this we construct
\[
\begin{CD}
S^{2np^{k+1}-1} \{p^{r+1}\} @> e >> \Ome P^{2np^{k+1}}(p^{r+1}) @ > \Ome \ol\dta_{k+1} >> \Ome F_n \\
@VV = V   @ VV = V   @ VV \Ome h V \\
S^{2np^{k+1}-1}\{p^{r+1}\} @ > e >> \Ome P^{2np^{k+1}}(p^{r+1}) @ > \Ome\ol\dta_k(np) >> \Ome F_{np}
\end{CD}
\]
where $k>0$, and 
\[
\begin{CD}
S^{2np-1}\{p^{r+1}\} @>>> \Ome P^{2np}(p^{r+1}) @ > \Ome \ol\dta_1 >> \Ome F_n \\
@ VVV  @ VVV  @VV \Ome h V \\
S^{2np-1} @>>> \Ome S^{2np} @>>> \Ome F_{np}
\end{CD}
\]
using diagram \eqref{eq4.3}.  Multiplying these together in order gives the result.\qed

\begin{lemma}
The map $\ol\phi_n\co P_n\to F_n$ can be chosen so that $P_n=P_{np} \vee Q_n$ and there is a homotopy commutative diagram:
\[
\begin{CD}
P_{np} \vee Q_n @> \ol\phi_n >> F_n \\
@ V p VV   @ VV h V \\
P_{np} @ > \ol\phi_{np} >> F_{np}
\end{CD}
\]
where $p$ is the projection.
\label{lem6.3}
\end{lemma}

{\bf Proof}\qua  Since $h\co H_*(F_n;\mathbb{Z}/p)\to H_* (F_{np};\mathbb{Z}/p)$ is onto the same holds for
\[
U(L^{(0)} (n)) \cong H_* (\Ome F_n;\mathbb{Z}/p)\to H_* (\Ome F_{np};\mathbb{Z}/p)\cong U(L^{(0)} (np))
\]
in fact the generators $x_{ip}\in L^{(0)} (n)$ satisfy
\[
h_*(x_{ip} (n))=u_ix_i(np)
\]
where $u_i$ is a $p$--local unit (see \fullref{prop5.3}).

Now given a basis $\{x_\alp, \bta^{(r)} x_\alp\}$ for $L^{(\infty)}(np)$, each $x_\alp$ is a Lie bracket in the $x_i(np)$ and this element consequently lifts to the same bracket in $x_{ip}(n)$.  Thus these liftings are linearly independent and can be chosen as part of a basis.  They are all in the image of the mod $p^r$ Hurewicz homomorphism.  Thus we have chosen generators for $L^{(\infty)} (n)$ which sort into those which are lifting of the generators for $L^{(\infty)}(np)$ and the others.  Realizing these via the Hurewicz homomorphism gives the maps $\ol\phi_n$.  By a change in basis for $Q_n$, we can assume that the map $p\co P_n\to P_{np}$ is trivial on $Q_n$.  

We now use \fullref{lem6.3} to construct the following ladder of fibrations:  
\begin{equation}
\begin{CD}
\Ome F_n @ > \ol\pat_n >>  \ol F_n @>>> P_n @ >\ol\phi_n >> F_n \\
@VVV    @VV \gam^s V @VVV @ VV h^s V \\
\Ome F_{np^s} @ >\ol\pat_{np^s} >> \ol F_{np^s} @>>> P_{np^s} @ > \ol\phi_{np^s} >> F_{np^s} 
\label{eq6.4}
\end{CD}
\end{equation}
and we use \fullref{lem6.2} to construct compatible equivalences:
\begin{equation}
\begin{CD}
S^{2n-1}\tms V_n @>>> \Ome F_n @> \pat_n >> \ol F_n \\
@VV \phi^s V  @VV \Ome h^s V  @VV \gam^s V \\
S^{2np^s-1} \tms V_{np^s} @>>> \Ome F_{np^s} @>>> \ol F_{np^s}.
\end{CD}
\label{eq6.5}
\end{equation}

Taking fibers vertically in \eqref{eq6.4}, we have a fibration sequence:
\[
\Ome \Dta_s\to K_s\to R_s\to \Dta_s
\]
where $R_s=(Q_n \vee \cdots \vee Q_{pn^s-1})\rtimes \Ome P_{np^s}$ is a wedge of mod $p^r$ Moore spaces and
\[
K_s=S^{2n-1}\tms \prod\limits_{1\le k<s} S^{2np^k-1} \{ p^{r+1}\} \tms\Ome S^{2np^s-1}.
\]
We use the splitting in \eqref{eq6.5} to construct a splitting
\[
K_s\to \Ome \Dta_s\to K_s.
\]
This completes the proof.\qed

\section{Decomposition of $\Ome F_{(p^s-1)}$}
\label{sec7}

Let $F_{(p^s-1)}$ be the $2np^s-1$ skeleton of $F$;  
since the mod $p^r$ homotopy classes
\[
\begin{CD}
P^{2ni-1} @>>> \Ome F_n
\end{CD}
\]
representing $x_i$ in homology factor through $F_{(p^s-1)}$ when $i<p^s$, the suspension map
\[
\begin{CD}
\sig\co  \wtit H_*(\Ome F_{(p^s-1)}; Z/p) @>>> H_* (F_{(p^s-1)}; k)
\end{CD}
\]
is onto.  By \fullref{propA3} of the Appendix, $H_*(\Ome F_{(p^s-1)}; \mathbb{Z}/p)$ is a tensor algebra on classes $x_i$ of dimension $2ni-1$ for $i<p^s$.  

\begin{theorem}
$\Ome F_{(p^s-1)}\simeq S^{2n-1} \tms \prod\limits_{1\le k<s} S^{2np^k-1} \{p^{r+1}\}\tms\Ome S^{2np^s-1} \tms \Ome Q_s$ where $Q_s\subset Q_{s+1}\subset \cdots \subset P_n$ is a wedge of mod $p^r$ Moore spaces for some choice of $P_n$ in \cite{CMN1}.
\label{thm7.1}
\end{theorem}

{\bf Proof}\qua  This proceeds along the same lines as \cite{CMN1}.  The maps 
\[
\ol\dta_k\co P^{2np^k}\to F
\]
factor through $F_{(p^s-1)}$ when $k<s$ as does the restriction
\[
\ol\dta_s\bigg|_{S^{2np^s-1}}\co S^{2np^s-1}\to F.
\]
We thus construct a map
\[
\Tht_s\co K_s=S^{2n-1}\tms \prod\limits_{1\le k < s} S^{2np^k-1} \{p^{r+1}\}\tms\Ome S^{2np^s-1} \to \Ome F_{(p^s-1)}
\]
as before.  The Bocksteins $\bta^{(i)}$ for $i<r$ are trivial and
$\Tht_s$ induces an isomorphism in the homology under the $r$th
Bockstein.  To see this observe that the Bockstein homology spectral
sequence of $H_*(\Ome F_{(p^s-1)};\mathbb{Z}/p)$ is a restriction of
the spectral sequence for $H_*(\Ome F;\mathbb{Z}/p)$ (see \fullref{secAppendix}).
$H_*(\Ome F_{(p^s-1)};\mathbb{Z}/p)=U(L^{(0)}_s)$ where $L^{(0)}_s$ is
a free Lie algebra generated by $x_i$ of dimension $2ni-1$ for $1\le i
< p^s$.

To complete the proof we construct Lie algebras $L^{(k)}_s\subset L_s^{(0)}$ which are compatible with the subalgebras $L^{(k)}$ of $L^{(0)}$.  First we examine the construction of $L^{(k)}$.  (See \fullref{sec2}).  Recall that we have short exact sequences:
\[
\begin{array}{ccccccccc}
0 &\lraw & L^{(1)} &\lraw & L^{(0)} &\lraw & \langle \tau_0 \rangle &\lraw & 0\\
0 &\lraw & L^{(k+1)} &\lraw & L^{(k)} &\stackrel{\tau_{k+1}}{\lraw} & \langle \tau_k, \sigma_k\rangle &\lraw & 0
\end{array}
\]
where $\pi_{k+1}$ is any map of Lie algebras such that
\begin{eqnarray*}
\pi_{k+1} (\sig_k) &=& \sig_k\\
\pi_{k+1} (\tau_k) &=& \tau_k
\end{eqnarray*}
where $\sig_k,\tau_k\in L^{(k)}\subset L^{(0)}$ are given by the formulas:
\begin{eqnarray*}
\tau_k &=& x_{p^k}\\
\sig_k &=& \frac{1}{2p} \sum\limits^{p^k-1}_{i=1} \binom{p^k}{i} [x_i,x_{p^k-i}].
\end{eqnarray*}
Nothing is said about the value of $\pi_{k+1}$ on the other generators.  We need to be more specific at this point.  Let us define the weight of an element in a free Lie algebra to be the minimal number of brackets in any term; in particular, $L^{(0)}$ is free on generators $x_i$, we define $\ome (x_i)=1, \ome [x,y]=\ome (x)+\ome (y)$ and $\ome(\Sig a_i)=\min \ome (a_i)$.  For an element $z\in L^{(k)}$ we define the weight of $z$ to be the weight considered as an element of $L^{(0)}$.  Thus $\ome (\tau_k)=1$ and $\ome(\sig_k)=2$.  We further specify the Lie algebra homomorphism $\pi_{k+1}$ by demanding that $\pi_{k+1}(z)=0$ if $\ome (z)>2$.  Now define $L^{(k)}_s = L_s\cap L^{(k)}$ for $k\le s$.  Since $\sig_k, \tau_k\in L_s^{(k)}$, we have short exact sequences:
\[
\begin{CD}
0 @>>> L^{(1)}_s @>>> L^{(0)}_s @>>> \langle \tau_0 \rangle @>>> 0 \\
  @VVV                @VVV                  @VVV            @| \\
0 @>>>  L^{(1)}  @>>> L^{(0)}   @>>> \langle \tau_0\rangle @>>> 0 \\
\\
0 @>>> L^{(k+1)}_s @>>> L^{(k)}_s @>>> \langle \sig_k, \tau_k\rangle @>>> 0 \\
@VVV  @VVV  @VVV  @| \\
0 @>>>  L^{(k+1)} @>>> L^{(k)} @>>> \langle \sig_k, \tau_k\rangle @>>> 0
\end{CD}
\]
for $k<s$.

Define $L_s^{(s+1)}$ to be the kernel of $L^{(s)}_s\to \langle \sig_s\rangle\to 0$.  In fact we have $L_s^{(s+1)}\subset L^{(\infty)}$.  To see this we need to show that $\pi_{r+1}(L_s^{(s+1)})=0$ for $r\ge s$.  The generators of $L_s^{(s+1)}$ of filtration 1 are of the form $x_i$ for $i<p^s$ and those of filtration 2 are of the form $[x_i,x_j]$ for $i,j<p^s$.  None of these have dimension $2np^s-1$, so $L_s^{(s+1)}$ lies in the kernel of the composition $L^{(s+1)}\to \langle \tau_s,\sig_s\rangle \to \langle \tau_s\rangle$.  Consequently $L_s^{(s+1)}\subset L^{(s+1)}$. Similarly for $r>s$ the generators of $L_s^{(s+1)}$ of weight 1 and 2 have dimensions $\le 4n (p^s-1)$ and consequently their images are 0 in $\langle \sig_r, \tau_r\rangle $ for $r>s$.

It follows that we may first choose a basis for $L_s^{(s+1)}$ and then choose a basis for $L^{(\infty)}$ containing these elements.  This completes the proof.\qed

\begin{corollary}
There is a commutative diagram of fibration sequences:
\[
\begin{CD}
\Ome F_n @>>> V_n @>* >> P_n @>>> F_n \\
@AAA                  @AAA        @AA A          @AAA  \\
\Ome F_{(p^{s+1}-1)} @>>> K_{s+1} @> * >> Q_{s+1} @>>> F_{(p^{s+1}-1)} \\
@ AAA   @ AAA   @ AAA  @ AAA \\
\Ome F_{(p^s-1)} @>>>  K_s @> * >> Q_s @ >>> F_{(p^s-1)}
\end{CD}
\]
where all the vertical maps are mod $p$ homology monomorphisms and each $Q_s$ is a wedge of mod $p^r$ Moore spaces.  
\label{cor7.2}
\end{corollary}

\section{On the sequence $F_{(p^s-1)}\to F_n \to F_{np^s}$}
\label{sec8}

Recall that there is a $p$--local fibration sequence
\[
S^{2n}_{(p^s-1)} \to S^{2n}_\infty \stackrel{H_{p^s}}{\lraw} S^{2np^s}_\infty.
\]
By analogy, we compare $F_{(p^s-1)}$ to the fiber $\Dta_s$ of the Hopf invariant $h^s$:
\[
\Dta_s\to F_n \stackrel{h^s}{\lraw} F_{np^s}
\]
\begin{theorem}
$\Ome \Dta_s \simeq \Ome F_{(p^s-1)} \tms \Ome X_s$ where $X_s$ is a wedge of mod $p^r$ Moore spaces.
\label{thm8.1}
\end{theorem}

\proof
For dimensional reasons, the inclusion $F_{(p^s-1)} \subset F_n$ lifts to $\Dta_s$. Consider the pull back diagram:
\[
\begin{CD}
\Ome Y_s @>>> \Ome F_{(p^s-1)} @>>> K_s @>>> Y_s @>>> F_{(p^s-1)} \\
@VVV   @VVV   @VV = V  @VVV   @VVV \\
\Ome R_s @>>>  \Ome\Dta_s @>>> K_s @> * >> R_s @>>> \Dta_s.
\end{CD}
\]
The map $K_s\to \Ome\Dta_s\to\Ome F_n$ constructed in \fullref{sec6} is obtained from maps:
\[
\begin{array}{lccccc}
\ol \dta_k\co P^{2np^k} &\lraw & F_n             &           \qquad\qquad          &  k<s \\
S^{2np^s-1}              &\lraw & P^{2np^s} &\stackrel{ \ol\dta_s}{\lraw}     F_n  & 
\end{array}
\]
and all of these maps factor through $F_{(p^s-1)}$ for dimensional reasons.  Consequently the retraction $K_s\to \Ome \Dta_s$ factors through $\Ome F_{(p^s-1)}$ as well and hence the map $K_s\to Y_s$ is null homotopic.  Consequently $\Ome Q_s\simeq \Ome Y_s$ and the map $Q_s\to F_{(p^s-1)}$ lifts to an equivalence $Q_s\simeq Y_s$.  Now the inclusion $Q_s\to P_n$ factors through $R_s$, so $R_s=Q_s\vee E_s$ and we have a diagram of fibrations:
\[
\begin{CD}
K_s @>>> Q_s @>>> F_{(p^s-1)} \\
@V=VV   @VVV   @VVV \\
K_s @>>> Q_s\vee E_s @>>> \Dta_s.
\end{CD}
\]
Now the fiber of $Q_s\to Q_s\vee E_s$ is $\Ome (E_s\rtimes \Ome Q_s)$ which is a retract of $\Ome (Q_s\vee E_s)=\Ome R_s$ and hence a retract of $\Ome\Dta_s$.  This completes the proof with $X_s=E_s\rtimes \Ome Q_s$.\endproof

\section{Factorization of $\ol\pat$}
\label{sec9}

In this section we will consider the components of
\[
\ol\pat_n\co \Ome^2 S^{2n+1}\lraw S^{2n-1} \tms V_n = S^{2n-1} \tms \prod\limits_{k\ge 1} S^{2np^k-1}\{p^{r+1}\}
\]
we will write $\pat_n(k)\co \Ome^2 S^{2n+1} \lraw S^{2np^k-1}\{p^{r+1}\}$ for the $k$th component, $k>0$ and $\ol\pat_n(0)$ for the projection $\Ome^2 S^{2n+1}\lraw S^{2n-1}$.

We begin by combining the diagram from \fullref{cor5.5}:
\[
\begin{CD}
\Ome^2 S^{2n+1}   @ > \pat_n >> \Ome F_n @ >>> \ol F_n \\
@V \Ome \wtit H_{p^s} VV  @VV \Ome h^s V   @VV h' V \\
\Ome^2 S^{2np^s+1} @> \pat_{np} >> \Ome F_{np^s} @>>> \ol F_{np^s} 
\end{CD}
\]
with the equivalences from \fullref{lem6.2}
\[
\begin{CD}
S^{2n-1}\tms V_n @> \Tht_n >> \Ome F_n @>>> \ol F_n \\
@VV \psi V  @VV \Ome h V   @VV h' V \\
S^{2np-1}\tms V_{np} @> \Tht_{np} >>  \Ome F_{np} @>>> \ol F_{np}
\end{CD}
\]
where $h'$ is defined in \fullref{lem6.3}, to get a homotopy commutative diagram:
\[
\begin{CD}
\Ome^2S^{2n+1}  @ > \ol\pat_n >> S^{2n-1} \tms V_n  \\
@V \Ome \wtit H_{p^s} VV    @VV \psi V \\
\Ome^2 S^{2np^s+1} @> \ol\pat_{np^s} >> S^{2np^s-1}\tms V_{np^s}.
\end{CD}
\]
From this we see that
\[
\ol\pat_n(k+1)=\ol\pat_{np^k}(1)\circ\Ome H_{p^k}\circ\gam_k, \qquad k\ge 0
\]
while
\begin{equation}
\begin{CD}
\Ome^2S^{2n+1} @> \ol\pat_n(1) >>    S^{2np-1}\{p^{r+1}\} \\
@V \Ome \wtit H_p VV   @ VV p V \\
\Ome^2 S^{2np+1} @> \ol\pat_{np} (0) >>     S^{2np-1}.
\end{CD}
\label{eq9.1}
\end{equation}
Putting these together gives:
\begin{equation}
\ol\pat_n =\ol\pat_n (0) +\sum\limits^\infty_{k=0} \ol\pat_{np^{k-1}}
(1)\circ 
\Ome H_{p^{k-1}} \circ\gam_k
\label{eq9.2}
\end{equation}
where $\gam_k\co \Ome^2S^{2n+1}\to\Ome^2S^{2n+1}$ is an equivalence.
Equation \eqref{eq9.2} proves\break
\mbox{\fullref{thm1.2}}.

\mbk\noi
{\bf Remark}\qua  If we had better control over the relationship between the Hopf invariants $h^s$ and the James--Hopf invariant, we could remove the factor $\gam_k$.

\begin{proposition}
Suppose $\ol\pat_n(1)\co \Ome^2S^{2n+1}\to S^{2np-1} \{p^{r+1}\}$ is null homotopic for each $n$.  Then

{\rm(a)}\qua $\ol\pat_n$ factors through $S^{2n-1}$

{\rm(b)}\qua The loops on the $p$th James--Hopf invariant: 
\[
\Ome H_p\co \Ome^2 S^{2n+1}\lraw \Ome^2 S^{2np+1}
\] 
is homotopic to a composition:
\[
\Ome^2 S^{2n+1} \lraw \Ome^2 P^{2np+1} \stackrel{\Ome^2 p}{\lraw} \Ome^2 S^{2np+1}
\]
{\rm(c)}\qua The fiber $B$ of $\ol\pat_{np} (0)\co \Ome^2 S^{2np+1}\to S^{2np-1}$ is a retract of $\Ome^2P^{2np+1}$

{\rm(d)}\qua $B\simeq BW_n$.
\label{prop9.3}
\end{proposition}

{\bf Proof}\qua  (a) follows directly from \eqref{eq9.2}.  To verify (b), note that by hypothesis
\[
\Ome H_p\co \Ome^2 S^{2n+1} \lraw \Ome^2 S^{2np+1}
\]
lifts so the fiber $B$ of $\ol\pat_{np}(0)$.  The following diagram provides a map from $B$ to $\Ome^2 P^{2np+1}$:
\[
\begin{CD}
B @>>> \Ome^2P^{2np+1} \\
@ VVV   @VVV \\
\Ome^2 S^{2np+1} @> 1 >>  \Ome^2 S^{2np+1} \\
@VV \ol\pat_{np}(0) V   @VVV  \\
S^{2np-1} @>>>  \Ome F_{np} 
\end{CD}
\]
(c) follows since $\Ome F_n\simeq S^{2np-1}\tms V_{np}\tms \Ome P_{np}$, so 
\[
\Ome^2 P^{2np+1} \simeq B\tms \Ome V_{np} \tms \Ome^2 P_{np}.
\]
To verify (d), note that $S^2BW_n$ is a retract of $S^2 \Ome^2 S^{2n+1}$ by \cite{G2}, so there is a map $S^2 BW_n \to P^{2np+1}$ which is non zero in homology, It's adjoint provides an equivalence:
$$
BW_n\to \Ome^2 P^{2np+1} \to B.\eqno{\qed}
$$

\section{Relationship with $BW_n$}
\label{sec10}

In this section we describe a factorization of $\ol\pat$ using the classifying space $BW_n$ for the double suspension.  Recall from \cite{G2}

\begin{proposition}
Suppose $\alp\co S^{2n}\to Y$ and $\bta\co S^{2n-1}Y\to S^{4n-1}$ are maps such that $\bta\circ S^{2n-1} \alp\co  S^{4n-1} \to S^{2n-1} Y \to S^{4n-1}$ is homotopic to the identity.  Then there is a fibration sequence:
\[
S^{2n-1} \stackrel{\ol\alp}{\lraw} \Ome Y \lraw B.
\]
This is natural in the following sense: Suppose we are also given 
\begin{eqnarray*}
&&\alp' \co  S^{2n} \lraw Z \\
&&\bta' \co  S^{2n-1} Z \lraw S^{4n-1} \\
&&f \co Y\lraw Z 
\end{eqnarray*}
such that $\bta\sim \bta' \circ S^{2n-1} f$ and $\alp'\sim f\alp $.

Then we have a homotopy commutative diagram:
\[
\begin{CD}
S^{2n-1} @>\ol\alp >> \Ome Y @>>> B \\
@V 1 VV   @VV \Ome f V   @VV \psi V \\
S^{2n-1} @>\ol\alp ' >> \Ome Z    @>>> B'
\end{CD}
\]
\label{prop10.1}
\end{proposition}

\begin{corollary}
There is a homotopy commutative diagram:
\[
\begin{CD}
S^{2n-1} @>>> \Ome^2 S^{2n+1} @>>> BW_n \\
@VVV   @VVV   @VVV \\
S^{2n-1} @>>> \Ome F_n @>>> B
\end{CD}
\]
\label{cor10.2}
\end{corollary}
Since $S^{2n-1}$ is a retract of $\Ome F_n$, $B\simeq V_n\tms \Ome P_n$ and we have:

\begin{corollary}
There is a homotopy commutative diagram:
\[
\begin{CD}
\Ome^2 S^{2n+1} @>>>  BW_n \\
@VV \pat V    @VVV \\
\Ome F_n @>>>  V_n.
\end{CD}
\]
\label{cor10.3}
\end{corollary}
In particular, $\ol\pat_n(k)\co \Ome^2S^{2n+1} \to S^{2np^k-1}\{ p^{r+1}\}$ factors through $BW_n$, proving \fullref{thm1.3}.

\appendix
\section{Appendix}\setobjecttype{App}
\label{secAppendix}

This appendix has two purposes.  We begin with a discussion of an observation that appears in the work of Cohen, Moore and Neisendorfer.  We repeat it here because it sheds light on their results and is easily extended to the case of the filtrations in \fullref{sec7}.  In addition, we prove a general result about the homology of certain loop spaces which does not seem to be in the literature.  
\mbk
We begin by looking at the homology Serre spectral sequence for the fibering
\[
\begin{CD}
\Ome F_n @>>> \Ome P^{2n+1} @>{\Ome\pi} >> \Ome S^{2n+1}.
\end{CD}
\]
This is a multiplicative spectral sequence with $$E^2_{p,q}\cong H_p
(\Ome S^{2n+1}; \mathbb{Z}/p)\otimes H_q(\Ome F_n; \mathbb{Z}/p).$$
Since all elements of $E^2_{p,0}$ are permanent cycles,
$E^2=E^\infty$. Restricting this fibration to $S^{2n}_{(k)}$ leads to
a fibration
\[
\begin{CD}
\Ome F_n @>>> E_k @>>> S^{2n}_{(k)}
\end{CD}
\]
and $H_* (E_k;\mathbb{Z}/p)\subset H_* (E_{k+1};\mathbb{Z}/p)\subset \cdots \subset H_* (\Ome P^{2n+1} ; \mathbb{Z}/p)$.  

Thus $H_*(\Ome P^{2n+1};\mathbb{Z}/p)$ is a filtered differential group using $\bta^{(r)}$.  It has a spectral sequence:
\[
E^2_{p,q} = H_p(\Ome S^{2n+1};\bta^{(r)})\otimes H_q (H_*(\Ome F_n;\mathbb{Z}/p), \bta^{(r)})
\]
converging to $H_{p+q}(\Ome P^{2n+1};\mathbb{Z}/p) = 0$ if $p+q>0$.  This gives a multiplicative spectral sequence:
\[
E^2_{p,q}=\mathbb{Z}_p[u_{2n}] \otimes H_q(H_*(\Ome F_n;\mathbb{Z}/p); \bta^{(r)}) \Raw 0.
\]
This spectral sequence has the same form as the Serre Spectral sequence for the homology of the fibration:
\[
\begin{CD}
\Ome^2S^{2n+1} @>>> P\Ome S^{2n+1} @>>> \Ome S^{2n+1}.
\end{CD}
\]
The differentials are forced by the multiplicative structure and consequently these spectral sequences are isomorphic.  From this we conclude 

\begin{proposition}
{\rm \cite{CMN1}}\qua 
$H_*(H_* (\Ome F_n; \mathbb{Z}/p);\bta^{(r)} )\cong H^* (\Ome^2 S^{2n+1};\mathbb{Z}/p)$ as algebras and co-algebras.
\label{propA1}
\end{proposition}

\begin{corollary}
$H_* (H_*(\Ome F_{(p^s-1)};\mathbb{Z}/p);\bta^{(r)})\cong H_*(\Ome S^{2n}_{(p^s-1)};\mathbb{Z}/p)$.
\label{corA2}
\end{corollary}
The later follows by restricting the spectral sequence to the first $p^s-1$ columns.

\begin{proposition}\label{propA3}
Suppose that $X$ is simply connected and $k$ is a field.  Suppose that the suspension map:
\[
\begin{CD}
\sig\co \wtit H_{r-1} (\Ome X; k) @>>> \wtit H_r (X; k)
\end{CD}
\]
is onto.  Choose a right inverse $\tau$ to $\sig$ and let $T_*\subset \wtit H_* (\Ome X; k)$ be the image of $\tau$.  Then $H_*(\Ome X; k)$ is the tensor algebra on $T_*$. 
\end{proposition}

{\bf Proof}\qua  We will use the Serre spectral sequence for the loop space fibration.  Note that it is a spectral sequence of $H_*(\Ome X; k)$ modules.  Observe that the suspension map can be factored 
\[
\begin{CD}
\wtit H_{r-1} (\Ome X; k)= E^{2}_{0,r-1} @>>> E^r_{0,r-1} @< d_r < \cong < 
E^r_{r,0} @>>> E^2_{r,0} = \wtit H^\cdot_r (X; k)
\end{CD}
\]
where the first homomorphism is onto and the last is $1-1$.  By hypothesis, the total composite is onto so the last homomorphism is an isomorphism.  In particular if $x\in E^2_{r,0}$, $d_i(x)=0$ for $i<r$ and $d_r(x)=[\tau (x)]\in E^\tau_{0,r-1}$ which is non zero.  Now define
\[
F(r)_p = \left\{
\begin{array}{cl}
H_p(X;k) & \mbox{\rm if $p=0$ or $p\ge r$}\\
0        & \mbox{\rm if $0<p<r$}.
\end{array}
\right. 
\]
We claim that there is a quotient module $Q(r)$ of $H_*(\Ome X; k)$ such that
\[
E^r_{p,q} \cong F(r)_p\otimes Q(r)_q.
\]
This is certainly true when $r=2$ with $Q(2)=H_* (\Ome X;k)$.  We proceed by induction.  

Because of this tensor product decomposition,
\[
d^r(f\otimes q)=d^r(f\otimes 1) (1\otimes q)=1\otimes \tau (f)q
\]
and hence $d^r_{p,q}=0$ unless $p=r$.  That is, all non zero differentials terminate on the fiber.  Furthermore
\[
\begin{CD}
d^r_{r,q}\co E^r_{r,q} @>>> E^r_{0,q+r-1}
\end{CD}
\]
is a monomorphism since otherwise the kernel would need to be in the image of some differential.  Consequently
\begin{eqnarray*}
\ker d_r &=& F(r+1)\otimes Q(r)\\
\lim d_r &=& F(r+1)\otimes T_*Q(r).
\end{eqnarray*}
Now define $Q(r+1)=Q(r)/T_*Q(r)$ which completes the induction.

We observe that the spectral sequence is completely determined by $\tau$.  We complete the proof by applying the comparison theorem.  Let $Y$ be a wedge of spheres together with an isomorphism:
\[
\begin{CD}
\phi \co H_* (Y;k) @>>> H_*(X;k).
\end{CD}
\]
We will define a map of spectral sequences:
\[
\begin{CD}
E^r_{p,q}(Y) @>>> E^r_{p,q} (X)
\end{CD}
\]
$\phi$ defines a isomorphism $F(\phi)\co F(r)(Y)\lraw F(r)(X)$.  Since $Y$ is a co-H space, there is a map $v\co Y\lraw S\Ome Y$; we define $\tau(Y)$ as the composite
\[
\begin{CD}
\wtit H_r (Y;k) @>v_* >> \wtit H_r (S\Ome Y;k)\cong\wtit H_{r-1}(\Ome Y;k).
\end{CD}
\]
Now define $\Tht\co H_*(\Ome Y;k)\lraw H_*(\Ome X;k)$ is the ring homomorphism determined by
\[
\Tht (\tau(Y)(x))=\tau(x)(\phi(x)).
\]
This defines a $H_*(\Ome Y;k)$ module homomorphism
\[
E^2_{p,q}(Y)\lraw E^2_{p,q}(X).
\]
However at each stage
\[
d^r(f\otimes q) = \left\{
\begin{array}{cl}
1\otimes \tau(f) q &\ \ r=p\\
0 &\ \ r\ne p.
\end{array}
\right.
\]
So this defines a map of spectral sequences.  Since they both converge to 0 and $\phi$ is an isomorphism, $\Tht$ is as well completing the proof.\qed

\bibliographystyle{gtart}
\bibliography{link}

\end{document}